\documentclass[11pt,reqno]{amsart}

\usepackage{marginnote,amssymb,mathrsfs,latexsym,verbatim,pdfsync,color,mathabx}

\newcommand{\bbC}{{\mathbb C}}

\newcommand{\bbN}{{\mathbb N}}
\newcommand{\bbR}{{\mathbb R}}

\newcommand{\bbZ}{{\mathbb Z}} 

\def\bbC{{\mathbb C}}

\def\bbN{{\mathbb N}}

\def\bbR{{\mathbb R}}

\def\bbZ{{\mathbb Z}}

\def\cC{{\mathcal C}}

\def\cF{{\mathcal F}}

\def\cI{{\mathcal I}}

\def\cO{{\mathcal O}}
\def\cP{{\mathcal P}}

\def\cS{{\mathcal S}}

\def\Re{\operatorname{Re}}

\def\la{\langle}
\def\ra{\rangle}

\def\eps{\varepsilon}
\def\vp{\varphi}
\def\ov{\overline}
\def\p{\partial}

\def\ms{\medskip}

\def\supp{\operatorname{supp}}

\def\CB{\color{black} }
\def\CR{\color{red} }

\def\Lps{\dot{L}^p_s}
\def\dotDelta{\dot{\Delta}}
\def\BMO{\operatorname{BMO}}

\newtheorem{thm}{Theorem}[section]
\newtheorem{prop}[thm]{Proposition}
\newtheorem{cor}[thm]{Corollary}
\newtheorem{lem}[thm]{Lemma}
\newtheorem{defn}[thm]{Definition}
\newtheorem{remark}[thm]{Remark}

\newtheorem{theorem}{Theorem}

\begin{document}

\title[Homogeneous
Sobolev spaces]{Fractional Laplacian, homogeneous Sobolev spaces and
  their realizations} 
\author[a. Monguzzi, M. M. Peloso, M. Salvatori]{alessandro Monguzzi, 
Marco M. Peloso, Maura Salvatori}\address{Dipartimento di Matematica e Applicazioni, Universit\`a degli Studi di
  Milano--Bicocca, Via R. Cozzi 55, 20126 Milano, Italy}
\email{{\tt alessandro.monguzzi@unimib.it}}
\address{Dipartimento di Matematica, Universit\`a degli Studi di
  Milano, Via C. Saldini 50, 20133 Milano, Italy}
\email{{\tt marco.peloso@unimi.it}}
\email{{\tt maura.salvatori@unimi.it}}
\thanks{{\em Math Subject Classification} 46E35, 42B35}
\thanks{All authors are  partially supported by the grant PRIN 2015
  {\em Real and Complex Manifolds: Geometry, Topology and Harmonic
    Analysis}, and are members of the Gruppo Nazionale per l'Analisi
  Matematica, la Probabilit\`a e le loro Applicazioni (GNAMPA) of the
  Istituto Nazionale di Alta Matematica (INdAM)} 

\begin{abstract}
We study the fractional Laplacian and the homogeneous Sobolev spaces
on $\bbR^d$, by considering two definitions that are  both  considered
classical.  We compare these different definitions, and show how they
are related by providing an explicit correspondence between these two
spaces, and show that they admit the same representation. Along the
way we also prove
 some properties of the fractional
Laplacian.
\end{abstract}
\maketitle

\vspace{-1cm}

\section{Introduction and statement of the main results}

The goal of this paper is to clarify a point that in our opinion has
been overlooked in the literature. 
Classically, the homogeneous Sobolev
spaces and the fractional Laplacian
are defined in two different ways.  
In one case, we consider the Laplacian $\Delta$ 
as densely defined, self-adjoint and
positive on $L^2(\bbR^d)$, and following Komatsu \cite{Komatsu}
it is possible to define the fractional powers $\Delta^{s/2}$ by means
of the spectral theorem, $s>0$.  
 Denoting by $\cS$  the space of Schwartz functions, one observes that
 for $p\in(1,\infty)$, 
 $\| \Delta^{s/2} \vp\|_{L^p(\bbR^n)}$ is a norm on $\cS$.  The  
 {\em homogeneous Sobolev space}  is the space
$\dot{W}^{s,p}$  defined as the closure of $\cS$ in such a 
norm.
The second definition is modelled on the classical Littlewood--Paley
decomposition
 of function spaces (see details below) and gives rise to
an operator, that we will denote by $\dot\Delta^{s/2}$, acting on spaces 
of tempered distributions {\em modulo polynomials}.  

Below, we describe the latter approach and show
how these two definitions are related to each other,  by showing that
they admit the same,
explicit realization.    We mention that
the analysis involved by the fractional Laplacian $\Delta^s$ has drawn
great interest in the latest years, beginning with the groundbreaking papers
\cite{Caffarelli,hitchhiker}, see also the recent papers
\cite{bourdaud2,bourdaud3,brasco-salort,Kwasnicki}.   We also mention that we were led to consider this
problem while working on spaces of entire functions of exponential
type whose fractional Laplacian is $L^p$ on the real line, \cite{MPS-Bernstein}.
\ms

Let $\cS$ denote the space of Schwartz functions,
$\cS'$ denote the space of tempered
distributions.  We denote by $\bbN_0$ the set of nonnegative integers
and
by  $\cP_k$ be the set of all polynomials of $d$-real
variables of degree $\le k$, with $k\in\bbN_0$.  We also set
$\cP_\infty=\bigcup_{k=0}^{+\infty} \cP_k$, the collection of all
polynomials. We also write  
$\ov\bbN=\bbN_0\bigcup\{+\infty\}$.

For $k\in\ov\bbN$,  we consider the equivalence relations  in $\cS'$ 
$$
f\sim_k g \Longleftrightarrow f-g\in \cP_k
$$
and we denote by $\cS'/\cP_k$ the space of tempered distribution modulo
polynomials in $\cP_k$, that is, the space of the equivalence classes with
respect the equivalence relation $\sim_k$.  If $u\in\cS'$ we denote by
$[u]_k$ its equivalent class in $\cS'/\cP_k$.  
For $k\in \bbN_0$, we denote by $\cS_k$
the space of all Schwartz functions $\vp$ such that 
$$
\int x^\gamma\vp(x)\, dx=0
$$
for any multi-index $\gamma\in\bbN^d$ with $|\gamma|\le k$ and set
$\cS_\infty =\bigcap_{k=0}^{+\infty} \cS_k$.
Each space $\cS_k$, $k\in \ov\bbN$ is a subspace
of $\cS$ that 
inherits the same topology of $\cS$. Moreover, it is easy to see that  
$$
\big(\cS_k\big)'=\cS'/\cP_k \,. 
$$

For $p\in[1,\infty]$ we denote by  $L^p(\bbR^d)$ (or, simply by
$L^p$) the standard Lebesgue space on $\bbR^d$ and we write
$$
\|f\|_{L^p}^p=\int |f(x)|^p\, dx \,,
$$
if $p<\infty$, and $\|f\|_{L^\infty}
=\operatorname{ess}\sup_{x\in\bbR^d} |f(x)|$.  

If $f\in L^1(\bbR^d)$ we define the Fourier transform 
$$
\hat f (\xi) = \int f(x) e^{-2\pi i x\xi}\, dx \,.
$$ 
The Fourier transform induces a surjective isomorphism 
$\cF:\cS\to\cS$ and
then 
it extends to a unitary operator $\cF: L^2\to
L^2$ and to a surjective  isomorphism $\cF:\cS'\to\cS'$.
For $u\in\cS'$, we shall also write $\widecheck u $ to denote $\cF^{-1} u$.

It is important  here to observe that if $\vp\in\cS_\infty$, then for
all $\tau\in\bbR$, 
$\cF(|\cdot|^\tau\widecheck\vp)\in \cS_\infty$ as well,  see e.g. 
\cite{Triebel,Grafakos}. 
\ms

\subsection{The 
fractional Laplacian $\dot\Delta^{s/2}$ and the homogeneous
Sobolev spaces $\Lps$}
We now recall the perhaps most standard 
definition of the fractional powers of the Laplacian and of the homogeneous
Sobolev spaces.  This material is classical and well known, see e.g.
\cite{Triebel,Grafakos}.   

We now have the following definition.
\begin{defn}\label{defn-homo-lap} {\rm
Let  $s>0$.  For $\vp\in \cS_\infty$, we define the 
{\em fractional Laplacian} of $\vp$ as
$$
\dotDelta^{s/2} \vp =\cF^{-1}(|\cdot|^{s}\widehat \vp) \,.
$$
}
\end{defn}

The next lemma follows from the arguments in \cite[Chapter
1]{Grafakos}, or \cite[5.1.2]{Triebel}.
\begin{lem}\label{lem1}
The operator
$$
\dotDelta^{s/2} :\cS_\infty\to \cS_\infty
$$
is a surjective isomorphism with inverse 
$$
\cI_s \vp= \cF^{-1} \big(|\xi|^{-s}\widehat\vp\big) \,.
$$ 
\end{lem}

Next, 
given $[u]\in \cS'/\cP_\infty$ we define another distribution in
$\cS'/\cP_\infty$ by setting for all $\vp\in\cS_\infty$ 
\begin{equation}\label{distribution}
\big\la\cF^{-1}(|\cdot|^s\widehat u),\vp\big\ra
:=\big\la u,\cF(|\cdot|^s\widecheck \vp)\big\ra,
\end{equation}
where, with an abuse of notation,\footnote{We warn the reader that, we
  shall denote with the same symbol $\la\cdot,\cdot\ra$ different
  pairings of duality, such as $\cS$ and $\cS'$, $\cS_\infty$ and
  $\cS'/\cP_\infty$, 
$L^p$ and $L^{p'}$,
  etc. The actual pairing of duality should be clear from the context
  and there should not be any confusion.}   we denote  by $\la\cdot,\cdot\ra$
also
the  pairing of duality between $\cS_\infty$ and $\cS'/\cP_\infty$.  
We remark that, clearly, the terms of both sides of 
 \eqref{distribution} are independent of the choice of the
 representative in $[u]$.
The fact that \eqref{distribution} defines 
a well-defined distribution follows at once from Lemma \ref{lem1}. 

Then, we extend the  definition of $\dotDelta^{s/2}$ to
$(\cS_\infty)'=\cS'/\cP_\infty$. 
\begin{defn}\label{defn-frac-lap-eq-classes} {\rm
We define the operator 
$$
\dotDelta^{s/2}: \cS'/\cP_\infty\to \cS'/\cP_\infty
$$
by setting, for any $[u]\in \cS'/\cP_\infty$, 
\begin{equation}\label{frac-lapl}
\dotDelta^{s/2}[u] =\cF^{-1}(|\cdot|^{s}\widehat u) \,,
\end{equation}
 and we call it 
the   fractional Laplacian (of
order $s$) of $[u]$.
}
\end{defn}

Since the term on the right hand side is independent of the choice of
the representative in $[u]$, we may simply write
$[\dotDelta^{s/2}u]$ in place of $[\dotDelta^{s/2}[u]]$.

\begin{remark}\label{remark-Laplacians}{\em 
We stress the fact that the terminology ``homogeneous Laplacian'' is
non-standard, and that such operator is typically simply called
fractional Laplacian. However,  one of the main goals of the present
paper is to consider the two standard different definitions of the
fractional Laplacian and of the corresponding homogeneous Sobolev
spaces and to study how they are related to each other.
 Thus, we need to
distinguish them, and we do so by introducing such terminology.
}
\end{remark}

For $1<p<\infty$, we define $\dot{L}^p$ as the space of all
elements in $\cS'/\cP_\infty$ such that in every equivalence class there is a
representative that belongs to $L^p$, that is,
\begin{equation}\label{def-hom-Lp}
\dot{L}^p = \Big\{ [u]\in \cS'/\cP_\infty:\, \|[u]\|_{\dot L^p}^p
:=\inf_{P\in\cP_\infty}  \|u+P\|_{L^p}^{p} <+\infty \Big\} \,.
\end{equation}
Clearly, if $[u]\in \dot{L}^p$, the
 representative of $[u]$ in $L^p$ is unique.

We observe that, since $\cS_\infty$ is dense in $L^p$ for all
$p\in[1,+\infty)$, see Lemma \ref{added-lem} below, we have that 
$\dot\Delta^{s/2}u \in \dot L^p$, $p\in(1,\infty)$, if (and only if) the
mapping
$$
\cS_\infty\ni \psi \mapsto \la u, \dot\Delta^{s/2} \psi\ra
$$
is bounded in the $L^{p'}$-norm.  In fact, in this case, there exists a unique $g\in
L^p$ such that
$ \la u, \dot\Delta^{s/2} \psi\ra = \la g, \psi\ra$ for all
$\psi\in\cS_\infty$ and we set
$$
\dot\Delta^{s/2} u = [g]\,. 
$$
\ms

Next, notice that  by definition $\dotDelta^{s/2} P =0$ for any
polynomial $P$.
The homogeneous Sobolev spaces $\Lps$ are defined as
follows.

\begin{defn}\label{defn-hom-sob-space} {\rm
 Let $s>0$  and let $1<p<\infty$. Then, the homogeneous
 Sobolev space $\Lps$ is defined as 
$$
\Lps = \big\{
[f]\in\cS'/\cP_\infty:\, \dotDelta^{s/2}f
\in \dot L^p
\big\}\,,
$$
 and we set
 \begin{equation}\label{hom-sob-norm}
  \|f\|^p_{\Lps}=\|\dotDelta^{s/2}f\|_{\dot{L}^p} \,.
 \end{equation} \ms
}
\end{defn}

Notice that, because of \eqref{def-hom-Lp}, equation \eqref{hom-sob-norm}
defines a norm on $\Lps$.  Thus,  the homogeneous
 Sobolev space $\Lps$ is a space of  equivalent classes in $\cS'/\cP_\infty$.\ms

\subsection{Fractional powers of $\Delta$ and the homogeneous Sobolev
  spaces ${\dot{W}}^{s,p}$} 

We now describe more in details the other definition of the homogeneous Sobolev
  spaces.
As me mentioned,
the operator $\Delta$ is densely defined, self-adjoint and
positive on $L^2(\bbR^d)$, and following Komatsu \cite{Komatsu}
it is possible to define the fractional powers $\Delta^{\alpha/2}$ for
all $\alpha\in\bbC$, with $\Re\alpha>0$.  Here we restrict our
attention to the case $\alpha=s>0$.  
For $\vp\in\cS$,
by the spectral theorem  we have that 
\begin{equation}\label{frac-lapl-def2}
\Delta^{s/2} \vp  =\cF^{-1}(|\cdot|^s \widehat \vp) \,.
\end{equation}

We now have a simple lemma to describe the basic properties of 
the action of $\Delta^{s/2}$ on $\cS$.

\begin{lem}\label{added-lem}
 Let $s>0$ and $p\in(1,\infty)$.  The following properties
  hold:
\begin{itemize}
\item[(i)] $\cS_\infty$ is dense in $L^p$;
\item[(ii)] for every $\vp\in\cS$, $\Delta^{s/2}\vp\in L^p$;\smallskip 
\item[(iii)] $\|\Delta^{s/2}\vp\|_{L^p}$ is a norm on $\cS$.
\end{itemize}
\end{lem}

\proof
(i) is   well-known and also easy to prove directly.  For (ii), for 
$\vp\in\cS$, write  $\Delta^{s/2}\vp = \Delta^{s/2}(I+\Delta)^{-M}
(I+\delta)^M\vp$.  If $M>s$ it is easy to check that
$\Delta^{s/2}(I+\Delta)^{-M}$ is a Fourier multiplier satisfying
Mihlin--H\"ormander condition (see \cite[Theorem
5.2.7]{Grafakos-cl}), hence bounded on $L^p$,
$p\in(1,\infty)$. Since $(I+\delta)^M\vp\in L^p$ for all $p$'s, (ii)
follows. In order to prove (iii), we only need to observe that if
$\vp\in\cS$ and 
$(I+\Delta)^M\vp=0$, then $\widehat\vp$ is a distribution supported at
the origin, hence $\vp=0$. 
\qed \ms

Then, we have the following  definition. 
\begin{defn}\label{dotWps-def}{\rm
For $1<p<\infty$ and $s>0$ we define ${\dot{W}}^{s,p}$, also
called the {\em homogeneous Sobolev space} as
the closure of $\cS$ 
with respect to  the norm  $\|\Delta^{s/2} \vp\|_{L^p}$.  More
precisely, given the equivalent relation on the space of $L^p$-Cauchy
sequences of Schwartz functions,  
$\{ \vp_n\}\sim \{ \psi_n\} $ if  $\Delta^{s/2}(\vp_n-\psi_n)\to 0$ as $n\to+\infty$, and
$[\{ \vp_n\}]$ denote the equivalent classes, then 
\begin{multline}
{\dot{W}}^{s,p}
= \Big\{ \big[\{\vp_n\}\big]: \ \{
\vp_n\}\subseteq\cS,\, \{ \Delta^{s/2}\vp_n\} \text{\ is a Cauchy
sequence in\ } L^p, \\
\text{with } 
\big\| \big[\{\vp_n\}\big]\big\|_{{\dot{W}}^{s,p}} = \lim_{n\to+\infty}
\| \Delta^{s.2} \vp_n \|_{L^p} 
\Big\}\,.\label{dotWps-eq}
\end{multline}
}
\end{defn}

 If $\{ \Delta^{s/2}\vp_n\}$ is  a
Cauchy
sequences in $L^p$, and we denote by $f$ its limit, we set
$$
\|f\|_{{\dot{W}}^{s,p}} = \lim_{n\to+\infty} \| 
\Delta^{s/2}\vp_n\|_{L^p} \,. \ms
$$

We remark that it is equivalent to define ${\dot{W}}^{s,p}$ 
as the closure of the $\cC^\infty$ functions with compact support,
space that we denote by $\cC^\infty_c$, with respect to the same norm
$\|\Delta^{s/2} \vp\|_{L^p}$. 
\ms

Thus, we have described two different notions of the fractional
Laplacian, $\dot\Delta^{s/2}$ and $\Delta^{s/2}$, that in particular
are  defined and taking values 
on completely different spaces.  Both these
notions can be considered as ``classical''.  Moreover, we have
described two different scales $\Lps$ and $\dot{W}^{s,p}$ of
homogeneous Sobolev spaces, for $s>0$ and $p\in(1,\infty)$.

The main goals of this note are two.  The first one
is to describe the spaces ${\dot{W}}^{s,p}$
explicitly and obtain a realization of them as space of functions
(see below).
The  second one is to  show how
the homogeneous Sobolev spaces ${\dot{W}}^{s,p}$ and $\dot{L}^p_s$ are
related to each other, providing in fact an explicit one-to-one and
onto correspondence between them. 
\ms

Finally, observe that, if $\vp\in\cS$, then $\dotDelta^{s/2}\vp =
[\Delta^{s/2}\vp]_\infty$.
\ms

\subsection{The realization spaces $E^{s,p}$}

Following G. Bourdaud \cite{Bourdaud, bourdaud2, bourdaud3}, if $k\in\ov\bbN$ and
$\dot{X}$ is a given  subspace  of
$\cS'/\cP_k$ which is a Banach space, such that the natural inclusion of
$\dot{X}$ into $\cS'/\cP_k$ is continuous,
we call {\em realization} of $\dot{X}$ a subspace $E$ of $\cS'$ such
that there exists a bijective linear map
$$
R: \dot{X}\to E
$$ 
such that $\big[ R[u] \big] = [u]$ for every $[u]\in\dot{X}$. 
We endow $E$ of the norm given by $\| R[u]\|_E = \| [u] \|_{\dot{X}}$.   
Obviously, $E$ becomes a Banach space with such norm.

Thus, a realization of
$\Lps$  is a space $E$ of tempered distributions 
 such that for each $f\in  E$, the equivalent  class
modulo polynomials  $[f]$ of $f$ is in $\Lps$ and conversely, for each
$[f]\in\Lps$ there exists a unique $\tilde f\in[f]$ such that
$\tilde f\in E$, and such correspondence 
is an isometry.  Our first
main result is that the spaces $E^{s,p}$ defined below are a realization of $\dot{W}^{s,p}$,
for $p\in(1,+\infty)$, $s>0$. 
The spaces are also a realization of the spaces $\Lps$, as was
indicated 
already in \cite{Bourdaud, bourdaud2}, and here we give a proof of this
latter fact.  As a consequence, we obtain a precise correspondence between
the spaces  $\dot{W}^{s,p}$ and $\Lps$.  
 
In order to describe such realization spaces $E^{s,p}$ we
need to recall the definition of $\BMO$ and of the
homogeneous Lipschitz spaces 
$\dot{\Lambda}^s$, for $s>0$.   We begin with the latter one.
For $s>0$ we write $\lfloor s \rfloor$
to denote its integer part. 

\begin{defn}\label{diff-op-def}{\rm
The
difference operators $D^k_h$
of increment
$h\in\bbR^d\setminus\{0\}$ and of order  $k\in\bbN$
are defined as follows.
If $k=1$, we write
$D_h$ in place of $D_h^1$, is
$D_h f(x)=f(x+h)-f(x)$,  
 and
then inductively, for $k=2,3,\dots$
$$
D_h^k f(x)=D_h \big[ D_h^{k-1}f\big](x)\,.
$$
}
\end{defn}

For a non-negative integer $m$, we denote by $\cC^m$ the space of
continuously differentiable functions of order $m$.

We recall the following well-known facts, see \cite[6.3.1]{Grafakos}.
\begin{prop}\label{elem-prop-finite-diff}

{\rm (i)} For  $k=2,3,\dots$ we have the explicit expression 
$$
D_h^k f(x) =\sum_{j=0}^k (-1)^{k-j} \begin{pmatrix} k\\
  j \end{pmatrix} f(x+jh)\,.
$$

{\rm (ii)} If $f\in \cC(\bbR^d)$ is  such that
  $D_h^kf(x)=0$ for all $x,h\in\bbR^d$, $h\neq0$, then $f$ is a
  polynomial of degree at most $k-1$, and conversely,  $D_h^kf(x)=0$
  for all polynomials of of degree  $\le k-1$.
\end{prop}
\ms

\begin{defn}{\rm 
Let 
$\gamma>0$. We define the {\em homogeneous Lipschitz space} 
$\dot{\Lambda}^\gamma$ as 
$$
\dot{\Lambda}^\gamma = \Big\{ [f]\in \cC/\cP_{\lfloor \gamma
  \rfloor}\,  :\, 
\| f\|_{\dot{\Lambda}^\gamma} := \sup_{x,h\in\bbR,\,  h\neq0}
\frac{|D_h^{\lfloor \gamma \rfloor +1} f(x)|}{|h|^\gamma} <+\infty \Big\}\,.
$$
}
\end{defn}

We remark that, for all $\gamma>0$,
if   $[f]_{\lfloor
  \gamma\rfloor}\in \dot{\Lambda}^\gamma$, then  $f$ is a function of moderate growth, so
that $\dot{\Lambda}^\gamma \subseteq \cS'/\cP_{\lfloor \gamma
  \rfloor}$.  Moreover, 
$f$ is in $\cC^{\lfloor
  \gamma\rfloor}$.
For these and other properties of Lipschitz spaces, see
e.g. \cite{Krantz-Lipschitz-paper,Triebel,Grafakos}. \ms

Next, we introduce the Sobolev-$\BMO$ spaces.
\begin{defn}{\rm 
The space $\BMO$ is the space of locally integrable
functions, modulo constants such that
$$
\| f\|_{\BMO}:= \sup_{x\in\bbR^d} \sup_{r>0} \frac{1}{|B_r|} \int_{B_r(x)}
  |f(y)-f_{B_r}|\, dy <+\infty \,. 
$$
where $B_r(x)=B_r$ denotes the ball centered at $x$ of radius $r>0$,
$|B_r|$ its measure, and $f_{B_r}$ the average of $f$ over $B_r(x)$.

For $k=1,2,\dots$, we define the Sobolev-$\BMO$ space as
$$
S_k(\BMO) = \big\{ f\in\cS'/\cP_k:\, \p_x^\alpha f \in\BMO\ \text{for
} |\alpha|=k 
\big\} \,,
$$
and set
$$
\| f\|_{S_k(\BMO)} = \sum_{|\alpha|=k} \| \p^\alpha f\|_{\BMO} \,. 
$$
} 
\end{defn}

\begin{defn}{\rm 
Given $f\in\cS'$,  if there
exists a sequence $\{\vp_n\}\subseteq\cS$ such that
\begin{itemize}
\item[{\tiny $\bullet$}] $\vp_n\to f$ in $\cS'$,
\item[{\tiny $\bullet$}] $\{ \Delta^{s/2} \vp_n\}$ is a Cauchy
  sequence in $L^p(\bbR^d)$, 
\end{itemize}
then we say that $\Delta^{s/2} f$ is defined as  the $L^p$-limit of
$\{\Delta^{s/2} \vp_n\}$ and  we set $\| \Delta^{s/2} f \|_{L^p}=\lim_{n\to+\infty} \|
\Delta^{s/2} \vp_n\|_{L^p}$.
} 
\end{defn}

\begin{remark}\label{convergence-rem}{\rm
We observe that the definition is well-given since if
$\{\vp_n\},\{\psi_n\}\subseteq\cS$,  are two sequences such that
$\vp_n,\psi_n\to f$ in $\cS'$ and 
$\{\Delta^{s/2} \vp_n\}, \{\Delta^{s/2} \psi_n\}$ are both  Cauchy in 
$L^p$, then for every $\eta\in\cS_\infty$,
$$
\la \Delta^{s/2} \vp_n -\Delta^{s/2} \psi_n, \eta\ra = 
\la \vp_n -\psi_n, \Delta^{s/2} \eta\ra \to 0
\,
$$
as $n\to\infty$.  Since $\cS_\infty$ is dense in $L^p$,
$\{\Delta^{s/2} \vp_n\}, \{\Delta^{s/2} \psi_n\}$ have the same
$L^p$-limit.
}
\end{remark}
\ms

We now are ready to define the spaces $E^{s,p}$ that we will
show to be the realization spaces for ${\dot{W}}^{s,p}$.
For a sufficiently smooth function $f$, we denote by $P_{f;m;x_0}$ the
Taylor polynomial of $f$ of degree $m\in\bbN_0$ at $x_0\in\bbR^d$.

\begin{defn}\label{Esp-def}{\rm For $s>0$ and $p\in(1,+\infty)$, we define the
  spaces $E^{s,p}$ as follows.

{\rm  (i)} Let $0<s<\frac dp$, and let $p^*\in(1,\infty)$ given by
$ \frac{1}{p^*}=\frac1p -\frac sd$.
Then, we define
 $$
E^{s.p}=\big\{f\in L^{p^*}:\,
\| f\|_{E^{s,p}} := \|\Delta^{s/2} f\|_{L^p}<+\infty \big\} \,.
 $$

{\rm (ii)}  Let $s>d/p$, $s-d/p\not\in\bbN$ and let
 $m=\lfloor s-d/p \rfloor$. Then, we define
$$
E^{s.p}
=\Big\{ 
f\in\dot{\Lambda}^{s-\frac dp}:\, 
\ P_{f;m;0}=0\,,\  
 \| f\|_{E^{s,p}} := \|\Delta^{s/2} f\|_{L^p} <+\infty\, \Big\}\,.
$$

{\rm (iii)} Let $s-d/p\in\bbN$ and set $m=\lfloor s-d/p \rfloor$. 
Let $B$   be a fixed ball in $\bbR^d$. Then, if $m=0$,
$$
E^{s.p}
=\Big\{ 
f\in \BMO: 
f_B=0, \quad
 \, \| f\|_{E^{s,p}} :=\| \Delta^{s/2} f\|_{L^p}<+\infty\, \Big\}\,,
$$
while, if $m\ge1$, 
 \begin{multline*} 
E^{s.p}
=\Big\{ 
f\in S_m(\BMO) \cap\cC^{m-1}:
{\rm (i)}\,  \ P_{f;m-1;0}=0\,,\
  {\rm (ii)} \, (\p_x^\alpha  f )_B=0\, \quad\text{for }
  |\alpha|=m\,,\\
 {\rm (iii)} \, \| f\|_{E^{s,p}} :=\| \Delta^{s/2} f\|_{L^p}<+\infty\, \Big\}\,.
\end{multline*} 
}
\end{defn}

\subsection{Littlewood--Paley decomposition}
In order to state our second main result we need the Littlewood--Paley
decomposition of $\Lps$, for whose details 
see e.g. to  
 \cite[6.2.2]{Grafakos}. 
Let $\eta\in \cC^\infty_c(\bbR^d)$ such that $\supp\eta\subseteq
\{\xi:\, 1/2\le |\xi|\le 2\}$, identically $1$ on the annulus
$\{\xi:\, 1\le|\xi|\le3/2\}$ and such that
$$
\sum_{j\in\bbZ}\eta(2^{-j}\xi)=1 \,, \qquad\qquad\qquad
\xi\in\bbR^d\backslash \{0\} \,.
$$
For $j\in\bbZ$ we set $\eta_j = \eta(2^{-j}\cdot)$  and 
for  $f\in\cS'$  we define\footnote{Classically, these operators are
  denoted by $\Delta_j$.  Given the recurrent appearance of 
$\Delta$ in this paper, we decided to use the notation $M_j$ instead.}
$$
M_j f=\cF^{-1}(\eta_j\widehat f).
$$
We observe that, due to the support conditions of the $\eta_j$'s, 
for all $j,k\in\bbZ$ we
have that
\begin{itemize}
\item[{\tiny $\bullet$}] $M_j\le (M_{j-1}+M_j+M_{j+1})M_j\le 3M_j$;
\item[{\tiny $\bullet$}] $M_kM_j=0$ if $|k-j|>1$.
\end{itemize}

Then, $M_j f=0$ for any $f\in\cP_\infty$, therefore the operators
$M_j$ are well-defined on $\cS'/\cP_\infty$ and taking values in $\cS'$, as
it is immediate to check. Thus, for simplicity of notation, we write
$M_ju$ instead of $M_j [u]_\infty$, and
observe that 
\begin{equation}\label{u-dec-S'}
[u]_\infty  =\sum_{j\in\bbZ} [M_ju]_\infty\,,
\end{equation}
where the convergence is in 
$\cS'/\cP_\infty$.  
Then, we have the following characterization of the spaces
$\Lps$, $\dot\BMO:=\BMO/\cP_\infty$, and
$\dot\Lambda^\gamma$, respectively.

For $1<p<\infty$ and $s>0$, $[f]_\infty\in
\Lps$, if and only if  $\big(\sum_{j\in\bbZ} \big(2^{js}|M_j f|\big)^2\big)^{1/2}
\in L^p$ and 
\begin{equation}\label{LP-charact-Lps}
\bigg\| \bigg(\sum_{j\in\bbZ} \big(2^{js}|M_j f|\big)^2\bigg)^{1/2}
\bigg\|_{L^p} \approx \|f\|_{\dot{L}^p_s} \,.
\end{equation}

It holds that $[f]_\infty\in\dot\BMO$ if and only if 
$\big(\sum_{j\in\bbZ} |M_j f|^2\big)^{1/2}
\in L^\infty$ and in this case
\begin{equation}\label{LP-charact-BMO}
\sup_{x\in\bbR^d}  \bigg(\sum_{j\in\bbZ} |M_j f|^2\bigg)^{1/2}
\bigg\|_{L^p} \approx \|f\|_{\dot\BMO} \,.
\end{equation}

For $\gamma>0$, we have that $[f]_{\lfloor \gamma\rfloor}
\in \cS'/\cP_{\lfloor \gamma\rfloor}$, we have 
$[f]_{\lfloor \gamma\rfloor} \in
\dot\Lambda^\gamma$ 
if and only if  $\sup_{j\in\bbZ} 2^{j\gamma} \|M_j
f\|_\infty <+\infty$ 
and  we have
\begin{equation}\label{LP-charc-Lip}
\sup_{j\in\bbZ} 2^{j\gamma} \|M_j
f\|_\infty \approx \|f\|_{\dot\Lambda^\gamma} \,,
\end{equation}
see e.g. \cite[Theorem 6.3.6.]{Grafakos} 
\ms

\subsection{Statement of the main results}

Our first result describes explicitly the elements of
$\dot{W}^{s,p}$, $s>0$, $p\in(1,\infty)$.

\begin{theorem}\label{elements-Wspdot-charact}
Let $s>0$, $1<p<\infty$,
and let $\{\Delta^{s/2} \vp_n\}$ be a Cauchy sequence in $L^p$, with
$\vp_n\in\cS$. 
Then, the following properties hold.
\begin{itemize}
\item[(i)] If $0<s<d/p$, then there exists a unique $g\in L^{p^*}$
  such that $\vp_n\to g$ in $L^{p^*}$ as $n\to+\infty$, 
where $\frac{1}{p^*} = \frac1p   -\frac sd$.\smallskip
\item[(ii)] If $s-d/p=m\in\bbN_0$, then there exists
  a unique $g\in S_m(\BMO)$
  such that $[\vp_n]_m\to [g]_m$  in $ S_m(\BMO)$ as $n\to+\infty$.\smallskip
\item[(iii)] If $s>d/p$ and $s-d/p\not\in\bbN$, 
then there exists
  a unique $g\in \dot\Lambda^{s-d/p}$
  such that $[\vp_n]_m\to [g]_m$  in
$\dot\Lambda^{s-d/p}$ as $n\to+\infty$,  where $m=\lfloor s-d/p \rfloor$.
\end{itemize}
In particular, $\dot{W}^{s,p}$ is a space of functions if $s<d/p$,
while 
$\dot{W}^{s,p}\subseteq\cS'/\cP_m$, where $m=\lfloor s-d/p \rfloor$,
if $s\ge d/p$.
\end{theorem}

\begin{theorem}\label{main-thm-realization-Wsp-dot}
For $s>0$ and $p\in(1,\infty)$, the spaces $E^{s,p}$ are realization
spaces of ${\dot{W}}^{s,p}$.  More precisely, the space 
$\dot{W}^{s,p}$ can be identified with a subspace of $\cS'/\cP_m$,
where $m=\lfloor s-d/p \rfloor$, and for each $[u]_m\in
\dot{W}^{s,p}$ there exists a unique $f\in [u]_m\cap E^{s,p}$ and
such that 
$$ 
\|\Delta^{s/2} f\|_{L^p} = \| [u]_m\|_{{\dot{W}}^{s,p}} \,.
$$
\end{theorem}

Recall that, given $B=B(x_0,r)$ a ball in $\bbR^d$, and an locally integrable
function $f$, 
we denote by $f_B$ the
average of $f$ over $B$. 
\begin{theorem}\label{main-thm-realization-Lps-dot}
Let $s>0$, $1<p<\infty$.  Then the following properties hold.
\begin{itemize}
\item[(i)] If $0<s<d/p$, given $[u]_\infty\in \Lps$ there exists a
  unique $f\in [u]_\infty$ and such that $f\in E^{s,p}$, and it is
  given by
$$
f = \sum_{j\in\bbZ} M_j u \,,
$$
where the convergence is in
$L^{p^*}$. 
\item[(ii)] If $s-d/p=m\in\bbN_0$, 
let $B=B(x_0,r)$ be any fixed ball in $\bbR^d$.
Then, if $m=0$, given
  $[u]_\infty\in \Lps$ consider $f\in [u]_\infty$ 
  given by
\begin{equation}\label{the-series-2.0}
f= f_0+f_1 - (f_0+f_1)_B := \sum_{j\le0} \big( 
M_j u - M_ju(0) \big) +  \sum_{j\ge1} 
M_j u - (f_0+f_1)_B 
 \,,
\end{equation}
where the first series converges uniformly on compact subsets of
$\bbR^d$, while the second one in the $\BMO$-norm.
Then,  $f\in\BMO$ with $f_B=0$, hence $f\in E^{s,p}$.  

If $s-d/p=m\in\bbN$, $m\ge1$, given $[u]_\infty\in \Lps$, let
 \begin{equation}\label{the-series-2.m}
f= f_0+f_1 - (f_0+f_1)_B := \sum_{j\le0} \big( 
M_j u - P_{M_ju;m;0} \big) +  \sum_{j\ge1} 
M_j u - (f_0+f_1)_B 
 \,,
\end{equation} 
where the both series converge uniformly on compact subsets of
$\bbR^d$, together with all derivatives up to order $m-1$.
Then, $f\in E^{s,p}$, in particular,
$f\in S_m(\BMO)$, $P_{f;m-1;0}=0$, and $(\p^\alpha f)_B=0$ for
$|\alpha|=m$. 
\item[(iii)] If $s>d/p$ and $s-d/p\not\in\bbN$, set $m=\lfloor
  s-d/p\rfloor$.
Then 
 given $[u]_\infty\in \Lps$ there exists a
  unique $f\in [u]_\infty$ such that $f\in E^{s,p}$, and it is
  given by
\begin{equation}\label{the-series-3}
f = \sum_{j\in\bbZ} \big( 
M_j u -P_{M_j u; m;0} \big)  \,,
\end{equation}
where the convergence is uniformly on compact sets in
$\bbR^d$, together with all the derivatives up to order $m$, and in the
$\dot\Lambda^{s-d/p}$-norm.
\end{itemize}
\end{theorem}

Finally, we provide the explicit isometric correspondence between the
two homogeneous Sobolev spaces $\dot{W}^{s,p}$ and $\Lps$. 
\begin{theorem}\label{Cor}
Let $s>0$, $1<p<\infty$ and let $m=\lfloor s-d/p\rfloor$.  
Then, given any $[u]_\infty\in\Lps$ there
exists a unique $f\in [u]_\infty$ such that $f\in E^{s,p}$, so that
$f\in\dot{W}^{s,p}$, or
$[f]_m \in \dot{W}^{s,p}$, resp., if $0<s<d/p$ or $s\ge d/p$, resp.,
with equality of norms.

Conversely, given $f\in\dot{W}^{s,p}$   if $0<s<d/p$, or 
$[f]_m \in \dot{W}^{s,p}$, 
 if $s\ge d/p$, resp., let $f$ be its image in the realization space
 $E^{s,p}$,  
then
$[f]_\infty\in\Lps$, with equality of norms.
\end{theorem}
 \ms
 
\section{Preliminary facts}\label{preliminaries-sec}

\subsection{The fractional integral operator $\cI_s$}
For $0<s<d$, we consider the Riesz potential $\cI_s$ defined for $\vp\in\cS$ as
$$ 
\cI_s \vp =
\cF^{-1} ( |\xi|^{-s}\widehat \vp) = \omega_{s,d} \big( |\cdot|^{s-d}* \vp \big)\,,
$$
where $\omega_{s,d}=\Gamma\big((d-s)/2\big)/2^s\pi^{s/2}\Gamma(s/2)$.
Observe that for $s\in(0,d)$, $|\xi|^{-s}$ is locally integrable, so
that if $\vp\in\cS$, $|\xi|^{-s}\widehat \vp\in L^1$ and $\cI_s$ is
well-defined. 

For $p\in(0,\infty)$ we denote by $H^p(\bbR^d)$, or simply by $H^p$,
the classical Hardy space on $\bbR^d$.  Having fixed $\Phi\in\cS$ with
$\int\Phi=1$, then
\begin{equation}\label{Hardy-sp-def}
H^p(\bbR^d) 
=\big\{ f\in\cS':\, f^*(x):=\sup_{t>0} |f*\Phi_t(x)| \in L^p(\bbR^d)
\big\} \,,
\end{equation}
where
$$
\|f\|_{H^p}= \|f^*\|_{L^p}\,.
$$
We recall that the definition of $H^p$ is independent of the choice of
$\Phi$ and that, when
$p\in(1,\infty)$, $H^p$ coincides with $L^p$, with equivalence of
norms.  We also recall that 
$\cS_\infty$ is dense in $H^p$ for all $p\in(0,\infty)$, see
\cite[Ch.II,5.2]{Stein}.
For these and other
properties of the Hardy spaces see e.g. \cite{Stein} or
\cite{Grafakos}.\ms

For the operator  $\cI_s$ the following regularity result holds, 
see \cite{Stein}, 
\cite{Grafakos}, \cite{Krantz-fractional-Hardy-paper} and in
particular \cite{Adams} for (ii).
\begin{prop}\label{Riesz-pot-prop}
Let $s>0$. The operator  $\cI_s$, initially defined for $\vp\in\cS$, extends to a
bounded operator
\begin{itemize}
\item[(i)] if $0<p<d/s$, then $\cI_s :H^p\to H^{p^*}$, where
  $\frac{1}{p^*} = \frac1p
  - \frac sd$;\smallskip
\item[(ii)] if $p=d/s$, then $\cI_s :H^p\to \BMO$;\smallskip
\item[(iii)] if $p>d/s$, then $\cI_s :H^p\to \dot{\Lambda}^{s-d/p}$.
\end{itemize}
\end{prop}
\ms

We also recall that the Riesz transforms are the operators $R_j$,
$j=1,\dots,d$,
 initially defined
on Schwartz functions
$$
R_j f(x) = \cF^{-1} \Big( 
-i\frac{\xi_j}{|\xi|} \widehat f\Big) \,, 
$$
and that they satisfy the relation $\sum_{j=1}^d R_j^2 =I$.  Moreover,
the $R_j$ are bounded 
$$
R_j :H^p \to  H^p
$$
for all $p\in(0,\infty)$, and $j=1,\dots,d$.

\ms

\section{Proofs of the main results}\label{realization-sec}  

Our first task is to describe explicitly 
the elements of  ${\dot{W}}^{s,p}$, that is, 
the limits of $L^p$-Cauchy
sequences of the form $\{\Delta^{s/2} \vp_n\}$, with $\vp_n\in\cS$.  
We are going to use this simple lemma, of which we include the
proof for sake of completeness.
\begin{lem}\label{integral-estimate}
Let $1<p<\infty$, $d/p<\nu<1+d/p$.  Given $a\in\bbR^d\setminus\{0\}$, let
and let $g_a(x)= |a-x|^{\nu-d} - |x|^{\nu-d}$.  Then, $g_a\in L^{p'}$
and 
$$
\| g_a\|_{L^{p'}} \lesssim |a|^{\nu-d/p}\,.
$$
\end{lem}

\proof
We first observe that $g_a\in L^{p'}_{{\rm loc}}$ if and only if
$\nu>d/p$, and then that
$\| g_a\|_{L^{p'}} =C |a|^{\nu-d/p}$, where
\begin{align*}
C^{p'} 
& = \int_0^\infty \int_S \Big| \frac{1}{|a'-rx'|^{d-\nu}} - \frac{1}{r^{d-\nu}}
\Big|^{p'} \, d\sigma(x')\, r^{d-1}\, dr\\
& = \int_0^\infty \frac{r^{d-1}}{r^{(d-\nu)p'}} \int_S \Big| \frac{1}{|a'/r-x'|^{d-\nu}} - 1
\Big|^{p'} \, d\sigma(x')\, dr \\
& =: \int_0^\infty \frac{r^{d-1}}{r^{(d-\nu)p'}}  F(r)\, dr \,,
\end{align*}
where $S$ denotes the unit sphere in $\bbR^d$, $d\sigma$ the induced
surface measure, and $a'=a/|a|$.  Thus, since by rotation invariance,
$C$ is independent of $a'$, 
we only need to show that $C$
is finite if and only if $0<\nu<1+d/p$. Notice that $F(r)$ is finite 
when $\nu>d/p$, so we only need to check the integrability of 
$\frac{r^{d-1}}{r^{(d-\nu)p'}}  F(r)$ as $r\to 0^+, \infty$. As $r\to
0^+$, $F(r)= \cO(1)$ and $\int_0^\infty
\frac{r^{d-1}}{r^{(d-\nu)p'}}<\infty$ if and only if $d/p<\nu$, which
is the case. As $r\to \infty$, $F(r)\approx 1/s^{p'}$, so that
$$
\int_2^\infty \frac{r^{d-1}}{r^{(d-\nu)p'}}  F(r)\, dr 
\approx \int_0^\infty \frac{1}{r^{(d-\nu)p'-d+1+p'}}  \, dr 
$$
which is finite if and only if $\nu-d/p<1$.  The conclusion now
follows. 
\qed
\ms

\proof[Proof of Thm. \ref{elements-Wspdot-charact}]
For $0<s<d$, it is immediate to check that for $\vp\in\cS$ we have
that
$\cI_s \Delta^{s/2} \vp=\vp$.  Therefore, by
Prop. \ref{Riesz-pot-prop} (i) we have
\begin{equation} \label{emb-ineq-1}
\|\vp\|_{L^{p^*}}
 = \| \cI_s \Delta^{s/2} \vp \|_{L^{p^*}} \lesssim \| \Delta^{s/2}
\vp \|_{L^p} \,.
\end{equation}
Then,  let $\{\vp_n\}\subseteq\cS$ be such that $\{
\Delta^{s/2}\vp_n\}$ is a Cauchy sequence in $L^p$.  By the above
estimate \eqref{emb-ineq-1} it follows that 
$$
\|\vp_m - \vp_n\|_{L^{p^*}}
  \lesssim \| \Delta^{s/2}\vp_m -\Delta^{s/2}
\vp_n \|_{L^p} \,,
$$
so that there exists $f\in L^{p^*}$ such that $\vp_n\to f$ in
$L^{p^*}$ as $n\to+\infty$. 
By Remark \ref{convergence-rem}, this implies that 
$\Delta^{s/2}\vp_n \to \Delta^{s/2}f$ in $L^p$.  Moreover,
\begin{equation} \label{emb-ineq-2}
\|f\|_{L^{p^*}} =\lim_{n\to+\infty} \| \vp_n\|_{L^{p*}} 
  \lesssim \lim_{n\to+\infty} \| \Delta^{s/2}
 \vp_n\|_{L^p} = 
\| \Delta^{s/2}f\|_{L^p} \,.
\end{equation}
This proves (i).

In order to prove case (ii) with $s=d/p$, 
let $\{\vp_n\}\subseteq\cS$ be such that $\{
\Delta^{s/2}\vp_n\}$ is a Cauchy sequence in $L^p$. Since
$0<s=d/p<d$,  
again we have that $\cI_s \Delta^{s/2} \vp_n=\vp_n$.
 By Prop. \ref{Riesz-pot-prop} (ii) it follows that 
$$
\|\vp_m -\vp_n\|_{\BMO}
  \lesssim \| \Delta^{s/2}\vp_m -\Delta^{s/2}
\vp_n \|_{L^p} \,.
$$
Hence, there exists a unique $g\in \BMO$ such that $\vp_n\to g$ in
$\BMO$, and it is such\ that $\Delta^{s/2}g\in L^p$.  For, given any
$\psi\in\cS_\infty$ we have $\Delta^{s/2}\psi\in \cS_\infty \subseteq
H^1$ so that
\begin{align} 
\la g, \Delta^{s/2}\psi \ra
& = \lim_{n\to+\infty} \la \vp_n ,\Delta^{s/2}\psi \ra 
=  \lim_{n\to+\infty} \la \Delta^{s/2}\vp_n , \psi \ra
= \la G,  \psi \ra \,, \label{bdd-fnct-Lp}
\end{align}
where $G$ is the $L^p$-limit of $\{ \Delta^{s/2}\vp_n\}$.  This shows
that the mapping $\psi\to \la g, \Delta^{s/2}\psi \ra$ is bounded in
the $L^{p'}$-norm on a dense subset, 
and therefore $\Delta^{s/2}g$ defines an element of $L^p$, that is $G$.

Let now $s-d/p=m\in\bbN$ and let $\{\vp_n\}\subseteq\cS$ be such that $\{
\Delta^{s/2}\vp_n\}$ is a Cauchy sequence in $L^p$.  Using the
identity
\begin{equation}\label{id-riesz}
| \Delta^{(s-1)/2} \nabla \vp|^2 = \sum_{j=1}^d |R_j \Delta^{s/2}
\vp|^2 \,,
\end{equation}
which is valid for all $\vp\in\cS$, and the boundedness of the Riesz
transforms, we see that $\{ \Delta^{(s-1)/2} \p_j \vp_n\}$ is a Cauchy
sequence in $L^p$, for $j=1,\dots,d$.  If 
$s-d/p=1$, we apply the previous argument to each of the sequences 
$\{ \Delta^{(s-1)/2} \p_j \vp_n\}$, $j=1,\dots,d$.
If $s-d/p=2$, the argument is analogous but simpler, 
given the identity $\Delta^{s/2}
= \Delta^{(s-2)/2}\Delta$, 
and in the general
case we argue inductively.  

We obtain that, for all multiindices
$\alpha$, with $|\alpha|=m$, there exists $g_\alpha\in\BMO$ such that 
$\p_x^\alpha \vp_n\to g_\alpha$ in $\BMO$ and
$\Delta^{(s-m)/2}g_\alpha \in L^p$.  Using the fact that in $\cS'$
$$
\p_{x_j} g_\alpha = \p_{x_\ell} g_\beta
$$
if $\alpha+e_j=\beta+e_\ell$, where $|\alpha|=|\beta|=$, 
$e_j,e_\ell$ are the standard
basis vectors, $j,\ell\in\{1,\dots,d\}$, and induction again, it is
easy to see that there exists $g\in S_m(\BMO)$ such that $\p^\alpha
g=g_\alpha$, for $|\alpha|=m$.   This implies that $\vp_n\to g$ in $ S_m(\BMO)$,
hence in particular in $\cS'/\cP_m$.
Let $\psi\in\cS_\infty$, so that we have $\Delta^{s/2}\psi\in \cS_\infty \subseteq
H^1$, so that
\begin{align*} 
\la g, \Delta^{s/2}\p^\alpha \psi \ra
& = \lim_{n\to+\infty} \la \p^\alpha \vp_n ,\Delta^{s/2} \psi \ra 
=  \lim_{n\to+\infty} \la \Delta^{s/2}\vp_n , \p^\alpha \psi \ra
= \la G,  \p^\alpha \psi \ra \,, 
\end{align*}
Hence, the mapping $\eta\to \la g, \Delta^{s/2}\eta \ra$ is bounded in
the $L^{p'}$-norm on $\cS_\infty$, that is,
  $\Delta^{s/2}g$ defines an element of
$L^p$.  This proves (ii).
\ms

Let now $s>d/p$ and $s-d/p\not\in\bbN$.  
  We assume first that $d/p<s<1+d/p$.  For $\vp\in\cS$ and $x\neq
y$ we write
\begin{align*}
\vp (x) -\vp(y)
& = \int \frac{e^{ix\xi}-e^{iy\xi}}{|\xi|^s} \big( \Delta^{s/2}
\vp)\widehat{\,}(\xi)\, d\xi\\
&= \int H(x,y,t)  \Delta^{s/2}
\vp(t) \, dt
\,,
\end{align*}
where
$$
H(x,y,t)
= \cF^{-1} \Big( \frac{e^{ix(\cdot)}-e^{iy(\cdot)}}{|\cdot|^s} \Big) (t)
= \gamma_{s,d}  \big( |t-x|^{s-d} - |t-y|^{s-d} \big) 
\,. 
$$
By Lemma \ref{integral-estimate} it follows that $H(x,y,\cdot)\in L^{p'}$ and that
$$
\|H(x,y,\cdot)\|_{L^{p'}}
\lesssim |x-y|^{s-d/p} \,.
$$
Therefore, 
\begin{align*}
\|\vp\|_{\dot{\Lambda}^{s-d/p}}
& \lesssim \| \Delta^{s/2} \vp \|_{L^p} \,.
\end{align*}
Hence, given $\{\Delta^{s/2}\vp_n\}$ is a Cauchy sequence in $L^p$,
there exists a unique $g\in \dot{\Lambda}^{s-d/p}$ such that $\vp_n\to
g$, as $n\to+\infty$. We argue as before, using the fact that 
the dual space of $H^q$ is $\dot\Lambda^{d(1/q-1)}$ when $0<q<1$, see
e.g. \cite[Ch. III]{GarciaCuerva-deFrancia}.
Let $q$ be given by $1/q= 1-1/p+s/d$, so that $0<q<1$. Then 
$\dot{\Lambda}^{s-d/p} = (H^q)^*$, and
given any
$\psi\in\cS_\infty$ we have $\Delta^{s/2}\psi\in \cS_\infty \subseteq
H^q$, so that
\begin{align} 
\la g, \Delta^{s/2}\psi \ra
& = \lim_{n\to+\infty} \la \vp_n ,\Delta^{s/2}\psi \ra 
=  \lim_{n\to+\infty} \la \Delta^{s/2}\vp_n , \psi \ra
= \la G,  \psi \ra \,, \label{bdd-fnct-Lp}
\end{align}
where $G$ is the $L^p$-limit of $\{ \Delta^{s/2}\vp_n\}$.  This shows
that the mapping $\psi\to \la g, \Delta^{s/2}\psi \ra$ is bounded in
the $L^{p'}$-norm on a dense subset, 
and therefore $\Delta^{s/2}g$ defines an element of $L^p$, that is $G$.
 We conclude that
$\Delta^{s/2}g\in L^p$ as well.  Finally, if $s-d/p\not\in\bbN$ and
$s-d/p>1$, let $m=\lfloor s-d/p \rfloor$, and $|\alpha|=m$.  Since,
$\| \Delta^{(s-m)/2}\p^\alpha \vp_n\|^{L^p} \le 
\| R^\alpha \Delta^{s/2}\p^\alpha \vp_n\|^{L^p} \le \|
\Delta^{s/2}\p^\alpha \vp_n\|^{L^p}$, $\{ \Delta^{(s-m)/2}\p^\alpha
\vp_n\}$ is a Cauchy sequence in $L^p$.
Using the first part
we can find $g_\alpha\in
\dot{\Lambda}^{s-m-d/p}$ such that 
$\p^\alpha \vp_n \to g_\alpha$ in
$\dot{\Lambda}^{s-m-d/p}$, for all $\alpha$ with $|\alpha|=m$.  
Hence, arguing as in (ii), there exists a unique
$g\in\dot{\Lambda}^{s-d/p}$ such that $\p^\alpha g=g_\alpha$ for 
$|\alpha|=m$.  Hence, 
$\vp_n\to g$ in $\cS'/\cP_m$,  since $\p^\alpha \vp_n \to \p^\alpha g$  in
$\dot{\Lambda}^{s-m-d/p}$
 for 
$|\alpha|=m$, $\vp_n\to g$ in $\dot{\Lambda}^{s-d/p}$.  Finally, by \eqref{bdd-fnct-Lp} it is
now clear that $\Delta^{s/2} g\in L^p$, and (iii) follows as well.
\qed
\ms

We state the characterization of ${\dot{W}}^{s,p}$, that also provides
the precise embedding result.  It may be considered {\em folklore} by
some, but to the best of our knowledge, there exists no explicit
statement in the literature.  
From this, the proof of Thm. \ref{main-thm-realization-Wsp-dot}
is then obvious.

\begin{cor}\label{real-hom-sob} For $s>0$ and $p\in(1,+\infty)$, the
 following properties hold.

\noindent
{\rm  (i)} Let $0<s<\frac dp$, and let $p^*\in(1,\infty)$ given by
$\frac{1}{p^*}=\frac1p - \frac sd$.
Then, 
 $$
{\dot{W}}^{s,p} =\big\{f\in L^{p^*}:\,
\|\Delta^{s/2} f\|_{L^p} <+\infty \big\} \,, 
 $$ 
and
$$
\| f\|_{ L^{p^*}} \lesssim \| f\|_{{\dot{W}}^{s,p}}\,.
$$

\noindent
{\rm (ii)}  
Let $s-d/p\in\bbN$ and set $m=\lfloor s-d/p \rfloor$. 
Then, 
$$
{\dot{W}}^{s,p} =\big\{f \in S_m(\BMO):\,
\|\Delta^{s/2} f\|_{L^p} <+\infty \big\} \,,
 $$ 
and
$$
\| f\|_{ S_m(\BMO)} \lesssim \| f\|_{{\dot{W}}^{s,p}}\,.
$$

\noindent
{\rm (iii)} 
Let $s>d/p$, $s-d/p\not\in\bbN$ and let
 set $m=\lfloor s-d/p \rfloor$.
Then, 
$$
{\dot{W}}^{s,p} =\big\{f \in \dot{\Lambda}^{s-d/p}:\,
\|\Delta^{s/2} f\|_{L^p} <+\infty \big\} \,,
 $$ 
and
$$
\| f\|_{ \dot{\Lambda}^{s-d/p}} \lesssim \| f\|_{{\dot{W}}^{s,p}}\,.
$$

In particular, $\dot{W}^{s,p}$ is a space of functions if $s<d/p$,
while 
$\dot{W}^{s,p}\subseteq\cS'/\cP_m$, where $m=\lfloor s-d/p \rfloor$,
if $s\ge d/p$.
\end{cor}

We also have the following characterization.  We point out that if $m$
is a negative integer, $\cP_m=\emptyset$ and $\cS'/\cP_m $ is $\cS'$
itself, as well as $\cS_m=\cS$.
\begin{cor}\label{ale-charact} 
For $s>0$ and $p\in(1,+\infty)$, let $m=\lfloor s-d/p\rfloor$. Then,
\begin{multline*}
\dot{W}^{s,p}
= \big\{ [f]_m\in\cS'/\cP_m :\, 
\text{(i) 
there exists a sequence}\ \{\vp_n\}\subseteq \cS\ \text{such that }
[\vp_n]_m\to [f]_m \text{ in }\cS'/\cP_m; \\
 \text{(ii)  the sequence} \ \{\Delta^{s/2}\vp_n\}\ \text{is a
   Cauchy sequence in}\ L^p \big\}\,.
\end{multline*}
\end{cor}

\proof
We have shown that any element of $\dot{W}^{s,p}$ is an element of
$\cS'/\cP_m$ satisfying the conditions $(i)$ and $(ii)$ above.
Conversely, let $[f]_m\in\cS'/\cP_m$ satisfy the conditions $(i)$ and
$(ii)$.  We need to show that for any $f\in[f]_m$, 
$\Delta^{s/2}f$ is a well-defined element of $L^p$ -- in fact the
$L^p$-limit of the sequence $\{\Delta^{s/2}\vp_n\}$.

Let $g=\lim_{n\to+\infty} \Delta^{s/2} \vp_n$ and fix
 $\psi\in\cS_\infty$. Then $\Delta^{s/2}\psi
\in\cS_\infty\subseteq\cS_m$ and we have
\begin{align*}
\la g,\psi\ra
&= \lim_{n\to+\infty} \la \Delta^{s/2} \vp_n, \psi\ra =
\lim_{n\to+\infty} \la \vp_n, \Delta^{s/2} \psi\ra 
= \lim_{n\to+\infty} \la [\vp_n]_m, \Delta^{s/2} \psi\ra \\
& = \la [f]_m ,\Delta^{s/2} \psi\ra \,.
\end{align*}
Therefore, we may set $\Delta^{s/2}[f]_m =g\in L^p$. 
\qed \ms

We also have the following, but significant observation.
\begin{cor}\label{Delta-s-ammazza-polimoni}
Let $s>0$, $m=[s-\frac d2]$.  Then, for every $P\in\cP_m$,
$\Delta^{s/2}P=0$. 
\end{cor}

\proof
We show that there exists sequences of Schwartz functions $\{\vp_n\}$
such that $\{\Delta^{s/2}\vp_n\}$ is a Cauchy sequence in $L^2$, while
$\vp_n\to P$ in $\cS'$.  Suppose first $s>\frac d2$ and
$s-\frac d2\not\in\bbN$. Let $\psi\in C^\infty_c$ be nonnegative,
identically $1$ for $|x|\le1 $ and with support in
$\{|x|\le2\}$. Define $\psi^\eps(x)=\psi(\eps x)$, and set $\vp_n =
P\psi^{\frac1n}$.  Then, it is obvious that $\vp_n \to
P$ in $\cS'$, while if we assume momentarily that $P=x^\alpha$, with
$|\alpha|=k$, and write 
\begin{align*}
\|\Delta^{s/2} \vp_n\|_{L^2}^2
& =  \int  |\xi|^{2s} \big| \big( \widehat P*
\widehat{\psi (\cdot/n)}\big)(\xi) \big|^2 \, d\xi 
 =  \int  |\xi|^{2s} \big| \big( \widehat P*
\widehat{\psi (\cdot/n)}\big)(\xi) \big|^2 \, d\xi \\
&   = n^{2(d+k)} \int  |\xi|^{2s} \big| 
\p_\xi^\alpha \widehat\psi (n\xi) \big|^2 \, d\xi 
 =  n^{d+2(k-s)} \int  |x|^{2s} \big| 
\p^\alpha \widehat\psi (x) \big|^2 \, dx \,,
\end{align*}
which tends to $0$ as $n\to+\infty$ when $s-\frac d2>k$, i.e. when
$k\le [s-\frac d2]$.

Next, we assume that $s=\frac d2$.  Let $\{\vp_n\}$ be a
sequence in $\cS_\infty$ such that $\Delta^{s/2} \vp_n \to 0$ in
$L^2$. Then, since $\vp_n \in \cS_\infty$,
$\vp_n=\cI_{d/2}\Delta^{s/2}\vp_n$, by  Prop. \ref{Riesz-pot-prop}
(ii), $\vp_n\to 0$ in $\BMO$, hence to a constant.  Finally, if
$s-\frac d2=m\in\bbN$, and we assume the conclusion valid for $m-1$,
we proceed inductively, using identity \eqref{id-riesz}.  If 
$ \Delta^{s/2} \vp_n \to 0$ in
$L^2$, then $\Delta^{s-1/2} \p_j\vp_n \to  0$ in
$L^2$ and therefore $\p_j \vp_n\to 0$ in $\cP_{m-1}$. This proves the lemma. 
\qed
\ms

We now turn to the proof of Thm. \ref{main-thm-realization-Lps-dot}.

\proof[Proof of Thm. \ref{main-thm-realization-Lps-dot}]
(i) Let $[u]_\infty\in \Lps$ with $s<d/p$, and set 
$\vp_n = \sum_{|j|\le n} M_j u$.  Then, $\vp_n\in\cS_\infty$ and using
Cor. \ref{real-hom-sob} (1) and the equality $\dot{\Delta}^{s/2}\vp=
[\Delta^{s/2}\vp]_\infty$ for functions in $\cS$ we have 
that, for $n<m\in\bbN$, 
\begin{align*}
\| \vp_m-\vp_n\|_{L^{p^*}} 
& = \Big\| \sum_{n<|j|\le m} M_j u\Big\|_{L^{p^*}} \lesssim
\Big\| \dot{\Delta}^{s/2} \Big( \sum_{n<|j|\le m} M_j u \Big) \Big\|_{L^p} 
= \Big\| \sum_{n<|j|\le m} M_j u\Big\|_{\Lps} \\
& \approx \Big\| \Big( \sum_{k\in\bbZ} \big( 2^{ks} \big| M_k
\textstyle{ \sum_{n<|j|\le m} M_j u } \big| \big)^2 \Big)^{1/2}
\Big\|_{L^p} \\
& \lesssim \Big\| \Big( 
\sum_{n-2\le |k|\le m+2}  \big( 2^{ks} \big| M_k u \big| \big)^2 \Big)^{1/2}
\Big\|_{L^p} 
\,.
\end{align*}
Since $u\in\Lps$, the equivalence of norms in
\eqref{LP-charact-Lps} implies that $\{\vp_n\}$ is a Cauchy
sequence in $L^{p^*}$.  Therefore, 
it follows that $f=\sum_{j\in\bbZ} M_j u$, where the convergence is in
$L^{p^*}$.  
This proves (i)
\ms

We now prove (iii).   Consider the series on the right hand side of
\eqref{the-series-3}. We claim that the series converges uniformly on compact sets in
$\bbR^d$, together with the series of partial derivatives up to order
$m:=\lfloor s-d/p \rfloor$.  Assuming the claim, and letting $f$ denotes
its sum, we clearly have that $f\in\cC^m$ and $P_{f;m;0}=0$. Moreover, 
letting $\Phi_n =\sum_{|j|\le n} (M_j u-P_{M_ju;m;0})$, $\vp_n=
  \sum_{|j|\le n} M_j u$, and using Cor. \ref{real-hom-sob}  we have 
\begin{align*}
\| \Phi_n \|_{\dot{\Lambda}^{s-d/p}} 
& = \| \vp_n \|_{\dot{\Lambda}^{s-d/p}} \lesssim \| \vp_n
\|_{\dot{W}^{s,p}} = \| \Delta^{s/2} \vp_n \|_{L^p} \\ 
& \le \Big\| \Big(\sum_{|k|\le n+1} 2^{ks} |M_ju|^2 \Big)^{1/2}
\Big\|_{L^p} \,.
\end{align*}
Hence, the series on the right hand side of \eqref{the-series-3}
converges also in $\dot{\Lambda}^{s-d/p}$ and uniquely determines $f$.
Then, the conclusion 
follows modulo the claim.  

The following estimate holds true for 
$f\in\dot{\Lambda}^\gamma$, $\gamma$ not an integer,
see \cite[Cor. 3.4]{Krantz-Lipschitz-paper}, 
\begin{equation}\label{vera?}
\big|f(x)-P_{f;m;y}(x)\big| \lesssim |x-y|^{\gamma-m} \| \nabla^m
f\|_{\dot{\Lambda}^{\gamma-\lfloor \gamma \rfloor}} \,,
\end{equation}
where we set $\nabla^m f = \big(\p^\alpha f\big)_{|\alpha|=m} $.
Since $s-m>d/p$, by  Sobolev's embedding theorem, Cor. \ref{real-hom-sob} (iii),
if
$|x|\le r$ we
have 
\begin{align*}
| \Phi_\ell (x) - \Phi_n(x) |
& \lesssim_r \| \nabla^m \Phi_\ell - \nabla^m \Phi_n\|_{\dot\Lambda^{s-m-d/p}} \lesssim \|
\vp_\ell - \vp_n\|_{\dot{W}^{s,p}}  \\
& \approx \Big\| \Big( \sum_{k\ge1} \Big(2^{ks} \big|M_k \big(
\sum_{\ell\le |j|\le n} \Delta_ju \big) \big|\Big)^2 \Big)^{1/2}
\Big\|_{L^p} \\
& \lesssim \Big\| \Big( \sum_{\ell-1\le k\le n+1} \big(2^{ks}
|M_k u| \big)^2  \Big)^{1/2} \Big\|_{L^p} 
\end{align*}
This shows that the series converges uniformly on
compact sets.  The same argument now easily applies to series of partial
derivatives up to order $m$, since
$$
| \p^\alpha \vp_\ell (x) - \p^\alpha \vp_n(x) |
 \lesssim_r\| \nabla^{m-|\alpha|} \p^\alpha \vp_\ell 
- \nabla^{m-|\alpha|} \p^\alpha \vp_n\|_{\Lambda^{s-m-d/p}} 
\lesssim \|\vp_\ell - \vp_n\|_{\dot{W}^{s,p}}\,.
$$
Therefore, $\Phi_n \to f$, which is an element of $[u]_\infty$, and
this proves (iii).
\ms

In order to prove (ii), 
let $s-d/p=m=0$, and
consider the two series on the right hand side of
\eqref{the-series-2.0}, that is, 
$\sum_{j\le0} \big( 
M_j u -M_j(0) \big)$, and $\sum_{j\ge1} 
M_j u$.  We begin with the latter one.  
We observe that $\sum_{j\ge1} 
M_j u$ converges to an element of $\dot\BMO:=\BMO/\cP_\infty$, 
say to $f+P$, for some
$P\in\cP_\infty$.  
However,
its Fourier transform has support in $\{|\xi|\ge1\}$, hence it must be $P=0$, so
that $\sum_{j\ge1} 
M_j u=f_1$ is a locally integrable function.

Next we show that the former series converge uniformly on compact subsets to
a function $f_0$.  
For $|x|\le r$
we have
\begin{align*}
|M_j u(x) -M_ju(0) | 
& \lesssim_r  \sup_{|x|\le r} |  M_j\nabla  u (x)|
 \lesssim  2^{j(1+d/p)}\| M_j  u \|_{L^p} \\
& \lesssim 2^j \Big( 
\int \Big( \sum_{k\in\bbZ}  |M_k 2^{jd/p}  M_j u (x)|^2 \Big)^{p/2} \, dx \Big)^{1/p} \\
& \lesssim 2^j \Big( 
\int \Big( \sum_{|k-j|\le1 } \big( 2^k |M_k   u(x)|\big) ^2
\Big)^{p/2} \, dx \Big)^{1/p} \\
& \lesssim 2^j \| [u]_\infty\|_{\Lps} 
\,.
\end{align*}
Since $j\le0$, the conclusion follows.  Now set
$$
f = f_0+f_1 - (f_0+f_1)_B \,,
$$
and the case $m=0$ is proved.

Let now $s-d/p=m\ge1$ and consider the series on the right hand side
of \eqref{the-series-2.m}.
For $|x|\le r$,  arguing as before 
we have
\begin{align*}
|M_j u(x) -P_{M_ju;m;0}(x) | 
& \lesssim_r  \sup_{|x|\le r} |  M_j\nabla^{m+1}  u (x)|
 \lesssim  2^{j(m+1+d/p)}\| M_j  u \|_{L^p} \\
& \lesssim 2^j \Big( 
\int \Big( \sum_{k\in\bbZ}  |M_k 2^{jd/p}  M_j u (x)|^2 \Big)^{p/2} \, dx \Big)^{1/p} \\
& \lesssim 2^{j(m+1)} \Big( 
\int \Big( \sum_{|k-j|\le1 } \big( 2^k |M_k   u(x)|\big) ^2
\Big)^{p/2} \, dx \Big)^{1/p} \\
& \lesssim 2^{j(m+1)} \| [u]_\infty\|_{\Lps} 
\,.
\end{align*}
Thus, the series converges uniformly on compact subsets.
Clearly, the same argument applies to all derivatives $\p^\alpha_x$ with
$|\alpha|\le m$.  Then, the function $f_0\in\cC^m$.  On the other
hand, the series $\sum_{j\ge1} M_ju$ is such that, for $|\alpha|\le
m-1$, recalling that $j\ge1$, 
\begin{align*}
|\p^\alpha M_ju(x)|
& \le 2^{j(|\alpha|+d/p)} \| M_ju\|_{L^p} \lesssim  
\Big( 
\int \Big( \sum_{|j-k|\le1} 
|M_k (2^{j(|\alpha|+d/p)}  M_ju)|^2 \Big)^{p/2} \, dx \Big)^{1/p}  \\
& \lesssim  
\Big( 
\int \Big( \sum_{|j-k|\le1} 
|M_k (2^{j(|\alpha|+d/p)}  M_ju)|^2 \Big)^{p/2} \, dx \Big)^{1/p} \\
& \lesssim  
2^{-j} \Big( 
\int \Big( \sum_{k\ge0} 
 (2^{ks} |M_k  u|) ^2 \Big)^{p/2} \, dx \Big)^{1/p} \,.
\end{align*}
Thus, also the series $\sum_{j\ge1} M_ju$ converges uniformly on
compact subsets, together with its partial derivatives $\p^\alpha$,
with $|\alpha|\le m-1$.  Finally, if $|\alpha|=m$ we have
$$
\p^\alpha( f_0 + f_1)
= \sum_{j\le0} M_j (\p^\alpha u) - M_j (\p^\alpha u)(0) +
\sum_{j\ge1} M_j (\p^\alpha u) \,,
$$
and we can repeat the argument of the case $m=0$. 

Finally, we set $f_0+f_1=f_2$ and define
$$
f= f_2- P_{f_2; m-1;0}(x) -
\sum_{|\alpha|=m} (\p^\alpha f_2)_B
\,.
$$
It is now easy to see that $f\in E^{s,p}$.
\qed
\ms

The proof of Theorem \ref{Cor} is now obvious, and we are done.
\ms

\section{Final remarks, and comparison with the work of Bourdaud}

In \cite{Bourdaud} Bourdaud proved a version of
Thm. \ref{main-thm-realization-Lps-dot}, and our work is inspired by
his.  The descriptions of the realization spaces are similar, but not
identical. In particular  we explicitly refer to the $\BMO$
space. Moreover,  Bourdaud is mainly focused on the techniques of the
Besov spaces, and the proof in the case of Sobolev spaces are
different, and typically more involved.  However, we point out that
Bourdaud proved  that the spaces $E^{s,p}$ with $s-d/p\not\in\bbN_0$
are the {\em unique} realization of $\Lps$ that are dilation
invariant, while, if $s-d/p \in\bbN_0$ there does not exist any
dilation
invariant realization of $\Lps$.

We mention that we have restrict ourselves to the case
$p\in(1,\infty)$.  The case of the Hardy spaces, that is, $p\in(0,1]$
is also of great interest and we believe it is worth further investigation.

Finally, 
it would also be very interesting to study analogous properties for
the homogeneous versions of Sobolev and Besov spaces in the
sub-Riemannian setting, see the recent papers
\cite{PV,BPTV,BPV,BPV-GAHA}.  \ms

{\em Acknowledgements.}  We warmly thank G. Bourdaud for several
useful comments and remarks 
  that improved the presentation of this paper. 

\bibliography{PW-bib-2}

\providecommand{\bysame}{\leavevmode\hbox to3em{\hrulefill}\thinspace}
\providecommand{\MR}{\relax\ifhmode\unskip\space\fi MR }
\providecommand{\MRhref}[2]{%
  \href{http://www.ams.org/mathscinet-getitem?mr=#1}{#2}
}
\providecommand{\href}[2]{#2}
\begin{thebibliography}{GCRdF85}

\bibitem[Ada75]{Adams}
David~R. Adams, \emph{A note on {R}iesz potentials}, Duke Math. J. \textbf{42}
  (1975), no.~4, 765--778. \MR{458158}

\bibitem[Bou88]{Bourdaud}
G.~Bourdaud, \emph{R\'{e}alisations des espaces de {B}esov homog\`enes}, Ark.
  Mat. \textbf{26} (1988), no.~1, 41--54. \MR{948279}

\bibitem[Bou11]{bourdaud2}
G\'{e}rard Bourdaud, \emph{Realizations of homogeneous {S}obolev spaces},
  Complex Var. Elliptic Equ. \textbf{56} (2011), no.~10-11, 857--874.
  \MR{2838225}

\bibitem[Bou13]{bourdaud3}
\bysame, \emph{Realizations of homogeneous {B}esov and {L}izorkin-{T}riebel
  spaces}, Math. Nachr. \textbf{286} (2013), no.~5-6, 476--491. \MR{3048126}

\bibitem[BPTV19]{BPTV}
Tommaso Bruno, Marco~M. Peloso, Anita Tabacco, and Maria Vallarino,
  \emph{Sobolev spaces on {L}ie groups: embedding theorems and algebra
  properties}, J. Funct. Anal. \textbf{276} (2019), no.~10, 3014--3050.
  \MR{3944287}

\bibitem[BPV19a]{BPV}
Tommaso Bruno, Marco~M. Peloso, and Maria Vallarino, \emph{{Besov and
  Triebel--Lizorkin spaces on Lie groups}}, ArXiv e-prints (2019).

\bibitem[BPV19b]{BPV-GAHA}
\bysame, \emph{{Potential spaces on Lie groups}}, ArXiv e-prints (2019).

\bibitem[BS19]{brasco-salort}
Lorenzo Brasco and Ariel Salort, \emph{A note on homogeneous {S}obolev spaces
  of fractional order}, Ann. Mat. Pura Appl. (4) \textbf{198} (2019), no.~4,
  1295--1330. \MR{3987216}

\bibitem[CS07]{Caffarelli}
Luis Caffarelli and Luis Silvestre, \emph{An extension problem related to the
  fractional {L}aplacian}, Comm. Partial Differential Equations \textbf{32}
  (2007), no.~7-9, 1245--1260. \MR{2354493}

\bibitem[DNPV12]{hitchhiker}
Eleonora Di~Nezza, Giampiero Palatucci, and Enrico Valdinoci,
  \emph{Hitchhiker's guide to the fractional {S}obolev spaces}, Bull. Sci.
  Math. \textbf{136} (2012), no.~5, 521--573. \MR{2944369}

\bibitem[GCRdF85]{GarciaCuerva-deFrancia}
Jos\'{e} Garcia-Cuerva and Jos\'{e}~L. Rubio~de Francia, \emph{Weighted norm
  inequalities and related topics}, North-Holland Mathematics Studies, vol.
  116, North-Holland Publishing Co., Amsterdam, 1985, Notas de Matem\'{a}tica
  [Mathematical Notes], 104. \MR{807149}

\bibitem[Gra14a]{Grafakos-cl}
L.~Grafakos, \emph{Classical {F}ourier analysis}, third ed., Graduate Texts in
  Mathematics, vol. 249, Springer, New York, 2014. \MR{3243734}

\bibitem[Gra14b]{Grafakos}
\bysame, \emph{Modern {F}ourier analysis}, third ed., Graduate Texts in
  Mathematics, vol. 250, Springer, New York, 2014. \MR{3243741}

\bibitem[Kom66]{Komatsu}
H.~Komatsu, \emph{Fractional powers of operators}, Pacific J. Math. \textbf{19}
  (1966), 285--346. \MR{0201985}

\bibitem[Kra82]{Krantz-fractional-Hardy-paper}
Steven~G. Krantz, \emph{Fractional integration on {H}ardy spaces}, Studia Math.
  \textbf{73} (1982), no.~2, 87--94. \MR{667967}

\bibitem[Kra83]{Krantz-Lipschitz-paper}
\bysame, \emph{Lipschitz spaces, smoothness of functions, and approximation
  theory}, Exposition. Math. \textbf{1} (1983), no.~3, 193--260. \MR{782608}

\bibitem[Kwa17]{Kwasnicki}
Mateusz Kwa\'{s}nicki, \emph{Ten equivalent definitions of the fractional
  {L}aplace operator}, Fract. Calc. Appl. Anal. \textbf{20} (2017), no.~1,
  7--51. \MR{3613319}

\bibitem[MPS19]{MPS-Bernstein}
A.~{Monguzzi}, M.~M. {Peloso}, and M.~{Salvatori}, \emph{{Fractional
  Paley--Wiener and Bernstein spaces}}, ArXiv e-prints (2019).

\bibitem[PV18]{PV}
Marco~M. Peloso and Maria Vallarino, \emph{Sobolev algebras on nonunimodular
  {L}ie groups}, Calc. Var. Partial Differential Equations \textbf{57} (2018),
  no.~6, Art. 150, 34. \MR{3858833}

\bibitem[Ste93]{Stein}
E.~M. Stein, \emph{Harmonic analysis: real-variable methods, orthogonality, and
  oscillatory integrals}, Princeton Mathematical Series, vol.~43, Princeton
  University Press, Princeton, NJ, 1993, With the assistance of Timothy S.
  Murphy, Monographs in Harmonic Analysis, III. \MR{1232192}

\bibitem[Tri10]{Triebel}
Hans Triebel, \emph{Theory of function spaces}, Modern Birkh\"{a}user Classics,
  Birkh\"{a}user/Springer Basel AG, Basel, 2010, Reprint of 1983 edition
  [MR0730762]. \MR{3024598}

\end{thebibliography}
\bibliographystyle{amsalpha}

\end{document}